\documentclass[11pt, reqno]{amsart}

\usepackage{amsmath, comment, amsthm, amssymb}
\usepackage{amscd}
\usepackage{graphicx}
\usepackage{ragged2e}
\usepackage{setspace}
\usepackage{abstract}
\usepackage{color}
\usepackage[foot]{amsaddr}

\usepackage{enumitem}
\usepackage[percent]{overpic}

\usepackage[pdftex, linktocpage=true]{hyperref}

\DeclareFontFamily{U}{mathx}{\hyphenchar\font45}
\DeclareFontShape{U}{mathx}{m}{n}{
      <5> <6> <7> <8> <9> <10>
      <10.95> <12> <14.4> <17.28> <20.74> <24.88>
      mathx10
      }{}
\DeclareSymbolFont{mathx}{U}{mathx}{m}{n}
\DeclareFontSubstitution{U}{mathx}{m}{n}
\DeclareMathAccent{\widecheck} {0}{mathx}{"71}

\usepackage[format=plain,labelfont=bf,up,width=.99\textwidth]{caption}

\theoremstyle{plain}
\newtheorem{thm}{Theorem}[section]
\newtheorem{lemma}[thm]{Lemma}
\newtheorem{prop}[thm]{Proposition}

\theoremstyle{definition}
\newtheorem{defn}[thm]{Definition}

\theoremstyle{remark}
\newtheorem{remark}[thm]{Remark}

\theoremstyle{plain}
\newtheorem*{thmHOV}{{\bf Theorem HOV}}

\newcommand{\R}{\mathbb{R}}
\newcommand{\N}{\mathbb{N}}
\newcommand{\C}{\mathbb{C}}
\newcommand{\Z}{\mathbb{Z}}
\newcommand{\D}{\mathbb{D}}

\newcommand{\bigO}{\mathcal{O}}

\newcommand{\s}{\mathbb{S}}
\newcommand{\K}{\mathcal{K}}

\def\proof{\par\medskip\noindent {\bf Proof.\ \ }}

\def\qed{\hfill $\square$\\ }

\def \hvec #1#2{\left(#1,\\ #2\right)\!}

\def \diag #1#2#3#4#5#6#7#8{\begin{CD} #1 @>#5>> #2\\ @V#6VV  @VV#7V\\ #3 @>#8>> #4 \end{CD}}
\def\He{{H\'enon }}

\title[A new proof of a theorem of Hubbard--Oberste-Vorth]{A new proof of a theorem of\\ Hubbard--Oberste-Vorth}
\author[Remus Radu and Raluca Tanase]{Remus Radu$^1$ and Raluca Tanase$^1$}
\address{$^1$Institute for Mathematical Sciences, Stony Brook University, Stony Brook, NY 11794.}
\email{rradu@math.stonybrook.edu, rtanase@math.stonybrook.edu}

\date{\today}

\begin{document}
%=============================================================================
\maketitle

\vspace{.35cm}
\begin{abstract}
\noindent {\sc abstract.} We give a new proof of a theorem of Hubbard--Oberste-Vorth \cite{HOV2} for H{\'e}non  maps that are perturbations of a hyperbolic polynomial and recover the Julia set $J^{+}$ inside a polydisk as the image of the fixed point of a contracting operator. We also give different characterizations of the Julia sets $J$ and $J^{+}$ which prove useful for later applications.
\end{abstract}
\vspace{.35cm}

%=============================================================================
%\tableofcontents
%=============================================================================

\section{Introduction} \label{sec:Intro}

Fixed point theorems have found a lot of applications in dynamical systems in higher dimensions. They are used in proving the existence of the local stable and the local unstable manifold of a hyperbolic fixed point, or the existence of local foliations in the presence of a dominated splitting of the tangent bundle over an invariant set of a $\mathcal{C}^k$ self-map of a Riemannian manifold. In this article we give a description of the global structure of the Julia sets $J$ and $J^+$ of a dissipative hyperbolic \He map in $\C^2$ as the unique fixed point of a contracting operator in an appropriate function space. This provides an alternative proof of a well-known theorem of Hubbard and Oberste-Vorth \cite{HOV2}, which was one of the starting points (along with \cite{HOV1}, \cite{FM}, \cite{BS1}, \cite{BS2}, \cite{FS}, etc.) of more than two decades of research in dynamics in several complex variables. 
The proof that we give strengthens slightly the result of the theorem, and some of the tools developed here have found further applications to the study of \He maps with a semi-parabolic fixed point or cycle \cite{RT1} and their perturbations \cite{RT2}. 
 
A complex \He map $H_{p,a}:\C^2\rightarrow \C^2$ is defined by
$
H_{p,a} \hvec{x}{y} = \hvec{p(x)+ay}{ax},
$
where $p$ is a monic polynomial of degree $d\geq 2$. In this normalization the \He map has constant Jacobian equal to $-a^{2}$, but any other representation would work. The \He map is a biholomorphism whenever $a\neq 0$ with inverse $H_{p,a}^{-1} \hvec{x}{y} = \hvec{y}{x-p(y/a)}/a$.
From the point of view of dynamics, the interesting objects to study are the sets of points with bounded forward and respectively backward orbits 
under the iterations of the \He map. Define the invariant subsets as in \cite{HOV1}, \cite{BS1}, and \cite{FS}:
\begin{equation*}
K^{\pm}=\left\{ \hvec{x}{y}\in \C^2 : \left\| H_{p,a}^{\circ
n}\hvec{x}{y} \right\| \ \mbox{remains bounded as}\
n\rightarrow\pm\infty\right\},
\end{equation*}
as well as $K=K^{-} \cap K^{+}$. Then let $J^{\pm}=\partial K^{\pm}$ be the topological boundaries and let $J=J^{-}\cap J^{+}$. The sets $J$ and $J^{\pm}$ are  called the {\it Julia sets} of the \He map. Define the escaping sets $U^{\pm}=\C^{2}-K^{\pm}$. In this paper we will consider only dissipative maps $H_{p,a}$ (that is $|a|<1$). In this situation, it is known that $K^{-}$ has no interior and so $K^{-}=J^{-}$ \cite{BS1}, \cite{FM}. Understanding $J^+$ on the other hand is a non-trivial task. If the \He map is hyperbolic and dissipative then the interior of $K^{+}$ consists of the basins of attraction of finitely many attractive periodic points \cite{BS1}. Each basin of attraction is a Fatou-Bieberbach domain (a proper subset of $\C^{2}$, biholomorphic to $\C^{2}$). The common boundary of the basins is the set $J^{+}$. 

Hubbard and Oberste-Vorth \cite{HOV2} studied the structure of the Julia sets $J$, $J^{+}$ and $J^{-}$ for \He maps which are small perturbations of a hyperbolic polynomial $p$. Polynomials and \He maps have some fundamental differences: polynomials are not injective whereas \He maps are, polynomials and their rate of escape functions have finitely many critical points, on the other hand \He maps do not have any critical points in the usual sense, but their associated rate of escape functions have infinitely many critical points. Starting from the polynomial $p$, Hubbard and Oberste-Vorth 
create some objects that carry bijective dynamics (projective and inductive limits), and use those to describe the dynamics of the \He map on its Julia sets (see \cite[Theorem~1.4]{HOV2}). Their proof relies on telescopes for hyperbolic polynomials and crossed mappings. We will give a new proof of the theorem for the sets $J$ and $J^{+}$ in the language of a fixed point theorem. We will recover the set $J^{+}$ inside the bidisk $\D_{r}\times \D_{r}$ as the image of the unique fixed point of a contracting graph-transform operator in some function space $\mathcal{F}$, which we define in Section \ref{sec:spaceF}. We will complete the proof of the theorem in Section 
\ref{sec:conjugacy}, when we establish conjugacies between the \He map and certain model maps. We also obtain other new characterizations of the Julia sets $J$ and $J^+$. 
The construction resembles the proof of the Hadamard-Perron Theorem (see e.g. \cite{KH}). This approach has the advantage that it can be generalized to complex \He maps with a semi-parabolic fixed point \cite{RT1}, but the analysis in that case is much more complex (due to loss of hyperbolicity) and requires several delicate arguments. 

\smallskip
\noindent {\it Acknowledgements.} We thank John Hubbard for explaining us the details of \cite{HOV1} and \cite{HOV2}.
%=========================================================================
\section{Tools from one-dimensional dynamics}\label{sec:1D}

For a polynomial $p$ of degree $d\geq 2$, the {\it filled Julia set}  of $p$ is the set of points with bounded forward orbit
\[
    K_p = \{z\in \C : |p^{\circ n}(z)|\ \mbox{bounded as}\ n\rightarrow \infty \}.
\]
The set  $J_p=\partial K_p$ is the {\it Julia set} of $p$. As usual, $p^{\circ n} = p\circ p \circ \ldots \circ p$ denotes the $n$-th iterate of $p$.  
If $K_p$ is connected (or equivalently $J_{p}$ is connected) then there exists a unique analytic isomorphism 
\[
    \psi_{p}:\C-\overline{\D} \rightarrow \C-K_p
\]
such that $\psi_{p}(z^2)=p(\psi_{p}(z))$ and normalized so that $\psi_{p}(z)/z\rightarrow 1$ as $z\rightarrow \infty$. Furthermore, if $J_p$ is locally connected then the Riemann mapping $\psi_{p}$ extends to the boundary $\s^{1}$ and defines a continuous, surjective map $\gamma:\s^{1} \rightarrow J_p$. The boundary map $\gamma$ is called the Carath\'eodory loop. We refer to \cite{Mi} and \cite{DH} for more details.

An external ray $R_t$ is the image under the Riemann mapping $\psi_{p}$ of the straight line $\{re^{2\pi i t}, r>1\}$.  The Carath\'eodory loop is defined as $\gamma(t)=\lim\limits_{r\searrow 1}\psi_{p}( re^{2\pi i t})$ and we say that the ray $R_{t}$ lands at a point $\gamma(t)\in J_{p}$ if this limit exists. The external ray $R_0$ lands at the $\beta$-fixed point of $p$.  An equipotential for the polynomial $p$ is the image under the Riemann mapping $\psi_{p}$ of the circle $\{re^{2\pi i t}, t\in \R/\Z\}$ of radius $r>1$. 

A point $x$ is called a critical point of $p$ if $p'(x)=0$, in which case $c=p(x)$ is called a critical value. We say that $p$ is hyperbolic  if $p'$ is expanding on a neighborhood of the Julia set.

Throughout this paper we assume that $p$ is hyperbolic and has connected Julia set.  In this case, the filled Julia set $K_{p}$ is connected and locally connected, and none of the critical points of $p$ belong to the Julia set $J_{p}$ \cite{DH}. Moreover, all critical points of $p$ are attracted to attracting cycles, and the number of  attracting cycles is bounded above by $d-1$, by the Fatou-Shishikura inequality. For each attracting cycle, we consider a union $V_{i}$ of sufficiently small disks centered around the points of the cycle, such that $V_i$ is contained in the immediate basin of attraction and $p(V_{i})$ is relatively compact in $V_{i}$. Set $\Delta = \bigcup_{i=1}^k  V_i$, where $k$ is the number of attracting cycles. There exists a minimal iterate $n\geq 0$ such that $p^{-\circ n}(\Delta)$ contains all critical values of $p$. So $p^{-\circ (n+1)}(\Delta)$ belongs to the interior of the filled Julia set $K_p$ and contains all critical points of $p$. 

Consider the set
\begin{equation}\label{eq:U}
 U: = \C- \overline{p^{-\circ n}(\Delta)} - \{ z\in \C-K_{p} : |\psi_{p}^{-1}(z)|\geq R\}
\end{equation}
for some large $R>1$.

The set $U':=p^{-1}(U)\subset U$ is relatively compact in $U$, and $p:U'\rightarrow U$ is a degree $d$ covering map. Let $\mu$ be the Poincar\'e metric on $U$. The polynomial $p:U'\rightarrow U$ is strongly expanding with respect to the metric $\mu$. The construction of the sets $U$ and $U'$ is the same as in \cite{DH} and \cite{H}.

\begin{figure}[htb]
\begin{center}
\begin{overpic}[scale=1.1]{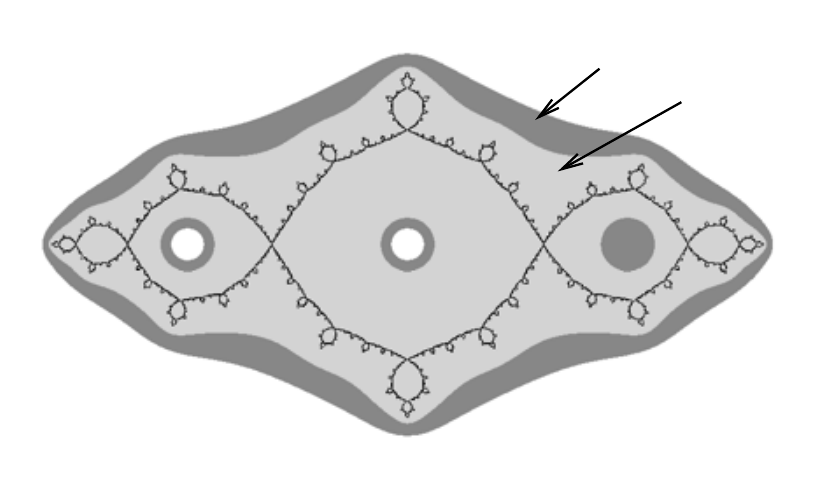}
 \put (81, 44.5) { {\large $U'$} }
 \put (71, 48.5) { {\large $U$} }
\end{overpic}
\end{center}
\vspace{-0.159cm}
\caption{A neighborhood $U$ of the Julia set of $p(z)=z^{2}-1$. The attracting cycle is $\{-1,0\}$  and $\Delta$ is a union of two small disks centered around the points of this cycle. The set $U$ (dark grey) is the complement of $\overline{\Delta}$ inside an equipotential of $p$, while $U'=p^{-1}(U)$ (light gray).}
\label{pic:nbd_U}
\end{figure}

 Choose $R$ as in Equation \ref{eq:U} and define the sequence of functions (equipotentials of the polynomial $p$) $\gamma_{n}:\R/\Z\rightarrow \C$ as follows:
\begin{equation}\label{eq:gamma}
    \gamma_{n+1}(t)=p^{-1}(\gamma_{n}(dt)):=\psi_{p}\left(R^{1/d^{n+1}}e^{2\pi i t}\right).
\end{equation}
With this notation, $\gamma_{-1}(\R/\Z)\subset\partial U$ and $\gamma_{0}(\R/\Z)\subset\partial U'$.

Since the Julia set $J_p$ is locally connected, the sequence of equipotentials $\gamma_n$ converges in the Poincar\'e metric of the set $U$ to the Carath\'eodory loop $\gamma$ of the polynomial $p$.

Let $\rho_{U}$ be the density function of the Poincar\'e metric on $U$, $\mu(z,dz)=\rho_U(z)|dz|$. The map $\rho_{U}$ is positive and $\mathcal{C}^{\infty}$-smooth on $U'$.  Since $U'$ is compactly contained in $U$, the Poincar\'e metric of $U$ is bounded below and above by the Euclidean metric on $U'$. If we let  $m=\inf\limits_{z\in U'}\rho_{U}(z)$ and $M=\sup\limits_{z\in U'}\rho_{U}(z)$ then
\begin{equation}\label{lemma: mM}
m|x-x'|\leq d_{U}(x,x')\leq M|x-x'|,
\end{equation}
for all $x, x' \in U'$. Consider now the constant $C := (\sup_{U'}|\rho_{U}'(z)|)/(\inf_{U'}\rho_{U}(z))$. The following lemmas will be useful later on.

\begin{lemma}\label{lemma:metric} Let $z$ be a point in  $U'$ and let $\delta$ be small enough so that $z-\delta$ is also a point in $U'$. Then $|\rho_{U}(z)-\rho_{U}(z-\delta)|\leq |\delta| C \rho_{U}(z)$.
\end{lemma}
The proof of the lemma is immediate and is left to the reader. 

\begin{lemma} \label{lemma: delta}
Let $z_{1}$ and $z_{2}$ be any two points in $U'$, and let $\delta$ be small enough
so that $z_{1}-\delta$ and $z_{2}-\delta$ are still in $U'$. Then
$
    d_{U}(z_{1}-\delta, z_{2}-\delta)\leq (1+|\delta| C) d_{U}(z_{1},z_{2}).
$
\end{lemma}
\proof
Let $\eta$ be a curve connecting $z_{1}$ and $z_{2}$, for which $\ell(\eta)=d_{U}(z_{1},z_{2})$. Then, if we translate $\eta$ by $\delta$, we get a curve (not necessarily length minimizing)
connecting $z_{1}-\delta$ to $z_{2}-\delta$. For small $\delta$, we can assume  that the new curve $\eta-\delta$ is still contained in $U'$. 
Its length is given by
\[
    \ell(\eta-\delta)=\int\limits_{\eta-\delta}\rho_{U}(z)|dz|=\int\limits_{\eta}\rho_{U}(z-\delta)|dz|.
\]
Using Lemma \ref{lemma:metric} we find that
\begin{eqnarray*}
\int\limits_{\eta}\rho_{U}(z-\delta)|dz| &\leq& \int\limits_{\eta}|\rho_{U}(z-\delta)-\rho_{U}(z)||dz| + \int\limits_{\eta}\rho_{U}(z)|dz| \\
&\leq& \int\limits_{\eta}|\delta| C \rho_{U}(z)|dz| + \int\limits_{\eta}\rho_{U}(z)|dz| = (1+|\delta| C)\ell(\eta).
\end{eqnarray*}
This shows that $d_{U}(z_{1}-\delta, z_{2}-\delta) \leq \ell(\eta-\delta) \leq (1+|\delta|C)\ell(\eta)=(1+|\delta|C)d_{U}(z_{1},z_{2})$.
\qed

%=============================================================================
\section{Construction of the neighborhood $V$}\label{sec:V}

Throughout this paper we will interchangeably use $H$ and $H_{p,a}$ to denote the \He map. 

By \cite{HOV1}, for $r$ sufficiently large, the space $\C^2$ can be divided
into three regions according to the dynamics of the \He map:
$ \D_r\times\D_r=\{(x,y)\in \C^2: |x|\leq r, |y|\leq r\},$
\[
W^+=\{(x,y):|x|\geq \max(|y|,r) \}\ \
\mbox{and}\ \ W^-=\{(x,y): |y|\geq \max(|x|, r)\}.
\]
\noindent The sets $J$ and $K$ are contained in the polydisk $\D_r\times\D_r$. The escaping sets $U^+$ and $U^-$ can be
described as union of backward iterates of $W^+$ and respectively forward iterates of $W^-$ under the \He map:
$ 
U^+=\bigcup_{k\geq 0} H^{-\circ k}(W^+) \ \ \mbox{and}\  \  U^-=\bigcup_{k\geq 0} H^{\circ k}(W^-).
$

Let $U'$ be the neighborhood of $J_{p}$ previously constructed.
Set $V:=U'\times \D_{r}$ for some $r>0$, chosen so that:
\begin{itemize}
  \item[(i)] $\overline{H(V)}$ does not intersect the horizontal boundary of $V$, that is $|ax|<r$ for any $x\in U'$.
  \item[(ii)]  $J\subset V$.
  One can choose for instance $r>3$ so that $J\subset \D_r\times\D_r$ as above. Notice that $J \cap \D_r\times\D_r = J \cap V$, by construction.
  \item[(iii)] All points in 
    $H(V)-\D_r\times \D_r$ belong to the escaping set $U^{+}$.  One can choose $R$ sufficiently large in Equation \ref{eq:U}
 so that the circle $\partial\D_r$ is contained in the set $U' $. By part (i), any point in $V$ that does not remain in $\D_r\times\D_r$ under forward iteration of $H$ belongs to the set $W^+$ which is contained in $U^+$.
\end{itemize}
Furthermore, suppose $|a|$ is small enough so that:
\begin{itemize}
  \item[(1)] $r|a|<\inf_{x\in U'}|p'(x)|$.
  \item[(2)] $r|a|<dist(\partial U', \partial U)$. 
  In other words, the $r|a|-$neighborhood of $U'$ is compactly contained in $U$.
\end{itemize}
The set $V$ is a neighborhood of the Julia set $J^+$ restricted to $\C\times\D_r$. 

\begin{figure}[htb]
\begin{center}
\begin{overpic}[scale=1]{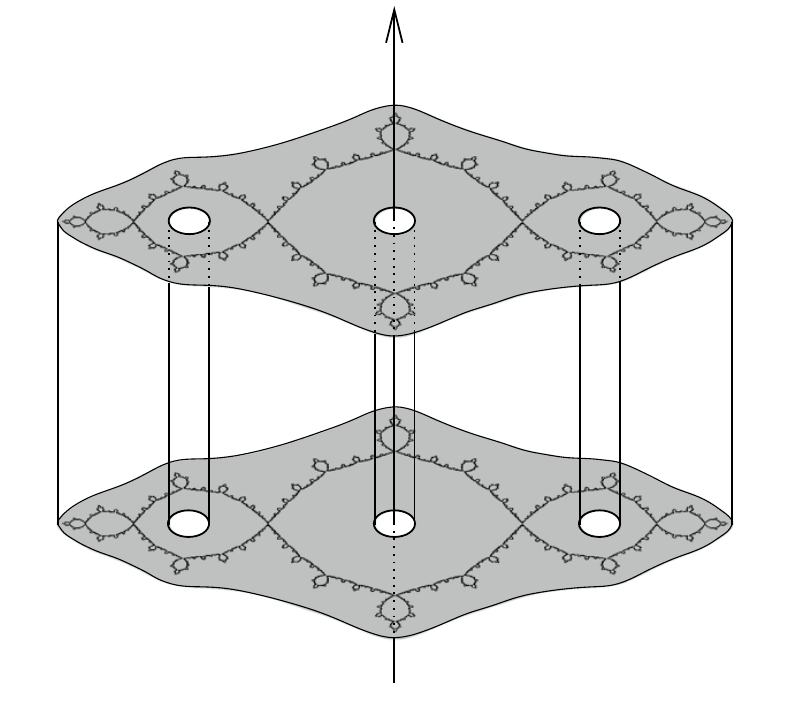}
 \put (70, 9) { {\large $U'$} }
 \put (93, 43) { {\large $\D_{r}$} }
 \put (52,81) {{\large $y$} }
\end{overpic}
\end{center}
\caption{A neighborhood $V=U'\times \D_{r}$ of $J^{+}\cap\{|y|<r\}$.}
\label{pic:V}
\end{figure}

\begin{lemma}\label{lemma:preim}
Let $(x,y)\in V$ and $(x',y')=H^{-1}(x,y)$. If $|y'|<r$ then $(x',y')\in V$.
\end{lemma}
\proof The point $\hvec{x'}{y'}$ belongs to $V$ iff $x'=y/a\in U'$ and $|y'|=\big{|}(x-p(y/a))/a\big{|}<r$.
By hypothesis we have that $\big{|}x-p(y/a)\big{|}<r|a|$. The point $x$ belongs to $U'$ and  $|a|$ is chosen small enough
so that the disk of radius $r|a|$ around $x$ is in $U$. It follows that $p(x')\in U$, hence $x'\in U'$. Therefore $(x',y')$ belongs to $V$.
\qed

\begin{prop}\label{prop: hyp-cones-E}
Let $(x,y), (x',y')$ be two points in $V$ with $H(x,y)=(x',y')$ and $(\xi,\eta)$ and $(\xi',\eta')$ two tangent vectors such that 
$DH_{(x,y)}(\xi,\eta)=(\xi',\eta')$. 
\begin{itemize}
\item[(a)] If $|\xi'|<|\eta'|$ then $|\xi|<|\eta|$.
\item[(b)] If $|\xi|>|\eta|$ then $|\xi'|>|\eta'|$.
\end{itemize}
\end{prop}
\proof
A direct computation gives $\xi'= p'(x) \xi+ a\eta$ and $\eta'=a\xi$.
\begin{itemize}
\item[(a)] If $|\xi'|<|\eta'|$ then $|p'(x)||\xi|-|a||\eta|<|\xi'|<|\eta'|=|a||\xi|$, so $|\xi|(|p'(x)|-|a|)<|a||\eta|$. The point $(x,y)$ belongs to $V$, so $x$ is bounded away from the critical points of $p$, in fact we have $|p'(x)|>r|a|$ where $r>2$.  Thus we get $|\xi|<|\eta|$.
\item[(b)] If $|\xi|>|\eta|$ then $|\xi'|>|p'(x)|\xi|-|a||\eta|>(|p'(x)|-|a|)|\xi|>|a||\xi|=|\eta'|$.
\qed
\end{itemize}

We define two invariant families of cones $\mathcal{C}^{h}_{(x,y)}$ and $\mathcal{C}^{v}_{(x,y)}$ in the tangent bundle of $V$,
\begin{eqnarray*}
\mathcal{C}^{h}_{(x,y)} &=& \{(\xi,\eta)\in T_{(x,y)}V :  |(x,\xi)|_{U}>|(y,\eta)|_{\D_{r}} \mbox{ and } |\xi|>|\eta| \}\\
\mathcal{C}^{v}_{(x,y)} &=& \{(\xi,\eta)\in T_{(x,y)}V : |(x,\xi)|_{U}<|(y,\eta)|_{\D_{r}} \mbox{ and } |\xi|<|\eta| \},
\end{eqnarray*}
where the lengths are measured with respect to the Poincar\'e metric on $U$ and $\D_{r}$, and with respect to the Euclidean metric. The cone invariance with respect to the Euclidean metric is shown in Proposition \ref{prop: hyp-cones-E}, whereas the invariance with respect to the Poincar\'e metrics has already been proven in \cite{HOV2}. We only use it to study vertical-like curves, so we will prove the part that we need at the end of Lemma \ref{lemma: two-preimages}.

\begin{defn} Let $\beta=\{(f(z),z),\ z\in \D_{r}\}\subset V$ be the graph of a holomorphic function $f:\D_{r}\rightarrow U'$. We say that $\beta$ is a vertical-like disk if  for all points $(x,y)$ on $\beta$, the tangent vectors to $\beta$ at $(x,y)$ belong to the vertical cone $\mathcal{C}^{v}_{(x,y)}$.
\end{defn}

\begin{lemma}\label{lemma: two-preimages}
If $\mathcal{\beta}$ is a vertical-like curve in $V$ then $H^{-1}(\mathcal{\beta})\cap V$ is the union of $d$ vertical-like curves.
\end{lemma}
\proof
By Lemma \ref{lemma:preim}, $H^{-1}(\beta)\cap V= H^{-1}(\beta)\cap \C\times\D_r$. 
Since the curve $\beta$ is vertical-like, it is the graph of a holomorphic function $f:\D_{r}\rightarrow U'$, hence $\beta=\{(f(z),z),\ z\in \D_{r}\}$.  
The function $f$ contracts Poincar\'e length and $|f'(z)|<1$. Then 
 \[
H^{-1}(\beta)=\left\{H^{-1}\hvec{f(z)}{z}=\hvec{z}{f(z)-p(z/a)}/a,\ z\in \D_{r}\right\}
\]
is an analytic curve whose horizontal foldings do not belong to the strip $\C\times \D_{r}$. Suppose there is a folding inside $\C\times \D_{r}$. Then, by Lemma \hyperref[lemma:preim]{\ref{lemma:preim}}, the folding point is actually inside $V$, hence its projection on the first coordinate $z/a$ belongs to $U'$ so it is bounded away from the critical points of $p$ (and the bound is independent of $a$). It follows that $p'(z/a)$ is bounded away from $0$, so $\frac{p'(z/a)}{a}$ gets arbitrarily large when $|a|$ is small enough whereas $f'(z)$ remains bounded, hence $f'(z)-\frac{p'(z/a)}{a}=0$ cannot have solutions inside $\D_{r}$.

Therefore the degree of the projection of $H^{-1}(\beta)$ on the second coordinate is constant in $\C\times \D_{r}$. It is easy to see that the degree of the projections is equal to the degree of the polynomial $p$, by looking at the number of intersections of $H^{-1}(\beta)$ with the $x$-axis. The curve $H(x\mbox{-axis})=\{(p(x),ax),x\in\C\}$ has $d$ connected components inside $V$, all horizontal-like. The curve $\beta$ is a vertical-like disk in $V$, hence $\beta$ intersects $H(x\mbox{-axis})$ in exactly $d$ points, which implies that $H^{-1}(\beta)$ intersects the $x\mbox{-axis}$ in $d$ points.

Thus $H^{-1}(\beta)\cap \C\times \D_r$ is a union of $d$ analytic curves $\beta_{i}$, $i=0,1,\ldots, d-1$,
which are all contained in $V$, by Lemma \hyperref[lemma:preim]{\ref{lemma:preim}}. 
The map $pr_{2}: \beta_{i}\rightarrow \D_{r},\ pr_{2}(x,y)=y$ is a covering map of degree one. By the Inverse Function Theorem,  $\beta_{i}$ is the graph of a holomorphic function $x=\phi(y)$ where $\phi:\D_{r}\rightarrow U'$. The map $\phi$ must also be injective, because $pr_{1}: \beta_{i}\rightarrow U',\ pr_{1}(x,y)=x$ is injective.  By the Schwarz-Pick lemma, $\phi:\D_{r}\rightarrow U'$ is weakly contracting in the Poincar\'e metrics of $\D_{r}$ and $U'$, hence strongly contracting if we endow $U'$ with the Poincar\'e metric of $U$. By Lemma \ref{prop: hyp-cones-E} we have $|\phi'(z)|<1$ for $z\in \D_{r}$. It follows that $\beta_{i}$ is vertical-like.
\qed

%=============================================================================
\section{A fixed point theorem}\label{sec:spaceF}

Consider the space of functions:
\begin{eqnarray*}
\mathcal{F} &=& \left\{    f:\s^{1}\times \D_{r} \rightarrow V  : 
			f(t,z)=(\varphi_t(z),z), \mbox{ where } f(t\times\D_{r})\ \mbox{is vertical-like,} \ \right. \\
   		   & &  \left.  \hspace{3.25cm} \varphi_t \mbox{ is analytic in } z \mbox{ and continuous in } t \right\}.
\end{eqnarray*}
We use the Kobayashi metric on $V$, which is simply the product of the Poincar\'e metric of U and the Poincar\'e metric of the vertical disk $\D_r$. On the function space $\mathcal{F}$ we consider the induced metric
\begin{eqnarray*}
 d(f, g) = \sup_{t\in \s^{1}}\sup_{z\in \D_{r}} d\left(pr_{1}(f(t,z)), pr_{1}(g(t,z))\right).
\end{eqnarray*}
The function space $\mathcal{F}$ is complete in the $d$-metric defined above.

\medskip
Let $\gamma_0$ be the equipotential of the polynomial $p$ (see Equation \ref{eq:gamma}) that defines the outer boundary of the set $U'$.  

\begin{defn} We  denote by $f_{0}:\s^{1}\times \D_{r}\rightarrow V$ the map $f_{0}(t,z)=(\gamma_{0}(t),z)$. The image of the map $f_0$ is a solid torus which represents the outer boundary of the set $V$.
\end{defn}

For any fixed $t\in \s^{1}$,  $f_{0}(dt\times\D_{r})$ is a vertical disk in $V$, so $H^{-1}\circ f_{0}(dt\times \D_{r})\cap V$ is a union of $d$ vertical-like disks, by Lemma \ref{lemma: two-preimages}. Let $C_{t}$ be the connected component of $H^{-1}\circ f_{0}(dt\times \D_{r})\cap V$ that crosses the $x$-axis at $\left(\gamma_{1}(t),0\right)$. Recall that $\gamma_1$ is the equipotential of the polynomial $p$ given by $\gamma_{1}(t)=p^{-1}(\gamma_{0}(dt))$, where the choice of the appropriate inverse branch of $p$ is made as in Equation \ref{eq:gamma}. Notice that $pr_{2}:C_{t}\rightarrow \D_{r}$, $pr_{2}(x,z)=z$ is a degree one covering map, hence $C_{t}$ is the graph of a holomorphic function $x=\varphi^1_{t}(z)$. This enables us to define a new function $f_{1}:\s^{1}\times \D_{r} \rightarrow V$ as $f_{1}(t,z)=(\varphi^1_{t}(z),z)$.  Notice that $f_{1}$ is homotopic to $f_{0}$ by construction since $\gamma_{1}$ and $\gamma_{0}$
are homotopic. Moreover, since $a$ is small, $f_{1}(\s^{1}\times \D_{r})$ and $f_{0}(\s^{1}\times \D_{r})$ are disjoint. Let $\tilde{\delta}=d(f_1, f_0)>0$. Notice that when $a$ is small $\tilde{\delta}$ is essentially the distance between $\partial U'$ and $\partial{U''}$ where $U''=p^{-1}(U')\Subset U'$.

\medskip
Let now $R_0:[0,1]\times\D_r\rightarrow V$, $R_0(0,z)=f_0(0,z)$, $R_0(1,z)=f_1(0,z)$ be a homotopy of vertical-like disks connecting $f_0(0\times \D_r)$ to $f_1(0\times \D_r)$, such that $R_0(s,0)$ is a point on the external ray of angle $0$ of the polynomial $p$ which connects $\gamma_0(0)$ to $\gamma_1(0)$.  As before, $H^{-1}(Im(R_0))\cap V$ has $d$ connected components. Denote by $R_1$ the component that contains  $f_1(0\times \D_r)$; $R_1$ is a collection of vertical-like disks 
that can be parametrized as graphs over the second coordinate, $R_1(s,z)=(\phi^1_s(z),z)$ for all $s\in[0,1]$. Inductively, we can construct a sequence of (approximative) external ray segments $R_n(s,z)=(\phi^n_s(z),z)$ by choosing the component of $H^{-1}(Im(R_{n-1}))\cap V$ that has the appropriate "matching end", i.e. for which $\phi^{n}_0(z)=\phi^{n-1}_1(z)$. The set $\mathcal{R}=\bigcup_{n\geq 0} R_n$ is our approximation for the external 3-D ray of angle $0$ for the \He map inside the set $V$.

\begin{defn}\label{def:F}
Consider now the subspace of functions $\mathcal{F'}\subset \mathcal{F}$, 
\begin{equation*}
\mathcal{F'} = \left\{ f_n:\s^{1}\times \D_{r} \rightarrow V : 
        f_0(t,z)=(\gamma_0(t),z),\ f_n(t,z)=F\circ f_{n-1}(t,z)\
        \mbox{for}\ n\geq 1
        \right\},
\end{equation*}
where the graph transform $F:\mathcal{F'}\rightarrow\mathcal{F'}$ is  defined as
\[
 F(f)=\tilde{f},
\] 
where $\tilde{f}$ and $f$ are homotopic and $ \tilde{f}\big{|}_{t\times\D_{r}}$ is the reparametrization  $\tilde{f}(t,z)=(\tilde{\varphi}_{t}(z),z)$ 
of the appropriate component of one of the $d$ vertical-like disk components of 
\[
H^{-1}\left(f(dt\times\D_{r}) \right)\cap V
\]
as a graph of a function over the second coordinate, via the Inverse Function Theorem.
\end{defn}

\begin{prop}
    The map $F:\mathcal{F'}\rightarrow\mathcal{F'}$ is well defined. 
\end{prop}
\proof Choose any function $f_n\in \mathcal{F'}$, $n\geq 1$.
The image of the map $f_n$ is a solid torus $T_1$ contained in the set $V$. The set \[
T_2=\{t\in\s^1, H^{-1}\left(f_n(dt\times\D_{r}) \right)\cap V\}
\]
is also a solid torus in $V$, which is mapped by the \He map to a solid torus wrapped around $d$ times inside $T_1$. In the $t$-coordinate, the \He map behaves like angle multiplication by a factor of $d$, while in the $z$ coordinate it acts like a strong contraction.  For each angle $t\in \s^1$, the set 
$\beta=f_n(dt\times \D_r)$ is a vertical-like disk in $V$. By Lemma \ref{lemma: two-preimages}, the set $L=H^{-1}\left( f_n(dt\times\D_{r}) \right)\cap V$ consists of $d$ vertical-like disks that we can label as $t+i/d$, for $i=0,1,\ldots, d-1$, and then parametrize as in Lemma \ref{lemma: two-preimages} as graphs over the second coordinate $ (\varphi^{n+1}_{t+i/d}(z),z)$. A choice of labelings that makes the map continuous with respect to $t$ is unique once we decide what the $0$-angle is for the new map. So we will call $f_{n+1}(0\times\D_r)$ the unique component of $H^{-1}\left( f_n(0\times\D_{r}) \right)\cap V$ that belongs the "external ray" $\mathcal{R}$. Then the map $f_{n+1}=F\circ f_n$ is simply defined as $f_{n+1}(t,z)=(\varphi^{n+1}_t(z),z)$ and is continuous with respect to $t$ and analytic with respect to $z$. 
\qed

\begin{thm}\label{thm: contraction}
The map $F:\mathcal{F'}\rightarrow\mathcal{F'}$ is a contraction in the metric defined on $\mathcal{F}$ and has an unique fixed point $f^*$.
\end{thm}
\proof
Consider any two functions  $f_n, f_k\in \mathcal{F'}$. We show that there exists a constant $\K<1$ such that, for any $t\in \s^{1}$:
\begin{equation*}
\sup\limits_{z\in \D_{r}} d_{U}\left( pr_{1}(F\circ f_{n}(t,z)), pr_{1}(F\circ f_{k}(t,z)) \right)\leq \K \sup\limits_{z\in \D_{r}}d_{U}\left( pr_{1}(f_{n}(dt,z)), pr_{1}(f_{k}(dt,z)) \right).
\end{equation*}

Recall that $f_{n}(dt\times \D_{r})$, $f_{k}(dt\times \D_{r})$,  $F\circ f_{n}(t\times \D_{r})$ and $F\circ f_{k}(t\times \D_{r})$ are vertical-like complex disks in $V=U'\times \D_{r}$, parametrized by the second coordinate, so there exists conformal maps $\psi_i, \varphi_{i}:\D_{r} \rightarrow U'$, for $i\in \{n,k\}$, such that $f_{i}(dt,z)=(\psi_{i}(z),z)$ and  $F\circ f_{i}(t,z)=(\varphi_{i}(z),z)$.

\begin{figure}[htb]
\begin{center}
\begin{overpic}[scale=0.95]{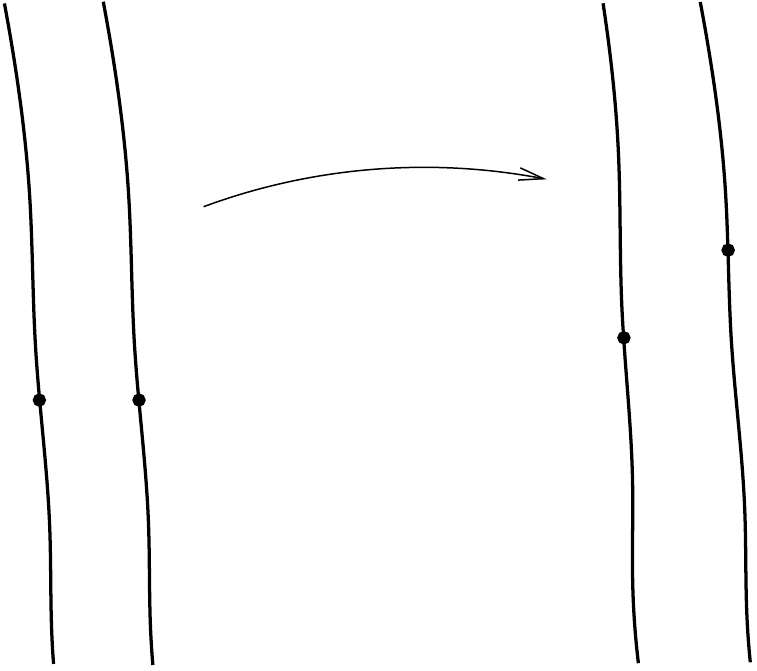}
 \put (45, 70) { { $H_{p,a}$} }
 \put (95, 55) { { $\left(p(x')+\delta,ax'\right)$} }
 \put (46, 43) { { $\left(p(x)+\delta,ax\right)$} }
 \put (17.5, 34.5) { { $\left(x',z\right)$} }
 \put (-12.5, 35) { { $\left(x,z\right)$} }
 \put (95, -6) { { $f_{k}(dt\times\D_{r})$} }
 \put (65, -6) { { $f_{n}(dt\times\D_{r})$} }
 \put (15, -6) { { $F\circ f_{k}(t\times\D_{r})$} }
 \put (-25, -6) { { $F\circ f_{n}(t\times\D_{r})$} }
\end{overpic}
\vspace{0.5cm}
\end{center}
\caption{Complex fibers $F\circ f_{n}$ and $F\circ f_{k}$ and their image under $H_{p,a}$.}
\label{pic:fibers}
\end{figure}
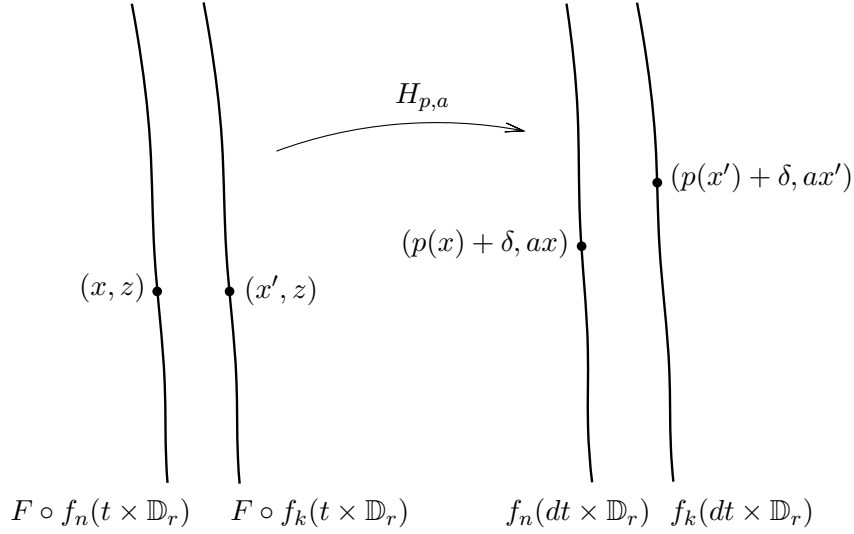

Let $z$ be any point in $\D_{r}$. Set $x= \varphi_{n}(z)$, $x'= \varphi_{k}(z)$, and $\delta = az$. Assume without loss of generality that $\delta<\tilde{\delta}$. With these notations we find that
\begin{eqnarray*}
H_{p,a} (x, z_{0}) &=& (p(x)+\delta, ax)=(\psi_n(ax),ax)\\
H_{p,a} ( x', z_{0})&=& (p(x')+\delta, ax')=(\psi_k(ax'),ax').
\end{eqnarray*}

The points $x$, $x'$, $p(x)+\delta$ and $p(x')+\delta$ all belong to $U'$.  Since $n,k\geq 1$ and $\delta<\tilde{\delta}$, the points $p(x)$ and $p(x')$ also belong to $U'$. The polynomial $p:U'\rightarrow U$ is strongly expanding with respect to the Poincar\'e metric of $U$, i.e. there exists a constant $\varepsilon<1$ (which depends only on the distance between  $\partial U$ and $\partial U'$) such that
\[
d_{U}(x,x')\leq \varepsilon d_{U}(p(x),p(x')).
\]
By Lemma  \ref{lemma: delta}, for small $\delta$, the following inequality holds:
\[
d_{U}(p(x), p(x'))\leq (1+|\delta| C) d_{U}(p(x)+\delta,p(x')+\delta).
\]
Thus we get
\begin{equation}\label{eq:KC}
d_{U}(x,x')\leq \varepsilon(1+|\delta|C) d_{U}(p(x)+\delta, p(x')+\delta).
\end{equation}

We now link the right hand side of  Equation \ref{eq:KC} with the distance between $f_{n}(dt\times\D_{r})$ and $f_{k}(dt\times\D_{r})$.
Notice that both fibers are vertical-like holomorphic disks, so the vertical distance between any two points of the fiber is bigger then their horizontal distance.  By the Schwarz-Pick lemma, the holomorphic map $\psi_n:\D_{r}\rightarrow U'$ is weakly contracting in the Poincar\'e metrics of $\D_{r}$ 
and $U'$, hence strongly contracting if we endow $U'$ with the Poincar\'e metric of $U$. It follows that 
\begin{eqnarray}\label{eq:1}
d_{U}(p(x)+\delta,p(x')+\delta)
&\leq& \sup\limits_{z\in \D_{r}}d_{U}(f_{n}(dt,z),f_{k}(dt,z)) + d_{U}(\psi_n(ax),\psi_n(ax')) \nonumber\\
&\leq& \sup\limits_{z\in \D_{r}}d_{U}(f_{n}(dt,z),f_{k}(dt,z)) + d_{\D_{r}}(ax,ax').
\end{eqnarray}

The set $H_{p,a}(\overline{V})$ does not intersect the vertical boundary of $V$, so $ax$ and $ax'$ belong to some disk $W$ compactly contained in $\D_{r}$. There exist constants $m_r$ and $M_r$ such that 
\[
    m_{r}|ax-ax'|\leq d_{\D_{r}}(ax,ax')\leq M_{r}|ax-ax'|.
\]

Following Equation \ref{lemma: mM}, a similar comparison holds if we put on the set $U'$ the Poincar\'e metric of $U$. Since $x,x'\in U'$ we get that
\[
    m|x-x'|\leq d_{U}(x,x')\leq M|x-x'|.
\]
Combining these two observations together with estimates \ref{eq:KC} and \ref{eq:1} we find that
\begin{equation*}
d_{U}(x,x')\leq \varepsilon (1+|\delta|C) \left( \sup\limits_{z\in \D_{r}}d_{U}(f_{n}(dt,z),f_{k}(dt,z))+ |a|\frac{M_{r}}{m}d_{U}(x,x') \right),
\end{equation*}
which yields
\begin{equation*}
d_{U}(x,x')\leq \K  \sup\limits_{z\in \D_{r}}d_{U}(f_{n}(dt,z),f_{k}(dt,z)),
\end{equation*}
where 
\[
\K := \frac{\varepsilon(1+|\delta|C) }{1-\varepsilon|a|(1+|\delta|C)\frac{M_{r}}{m} }.
\]
The constants $\varepsilon$, $C$, $m$ and $M_{r}$ are independent of $a$. The factor $\delta$ is small such that $|\delta|<|a|r$. Since $\varepsilon^{-1}>1$, there exists $a_{0}>0$ so that $1+|a_{0}|(1+rC)+|a_{0}|^{2}rC\frac{M_{r}}{m}<\varepsilon^{-1}$.  

Hence $\K<1$ for all $a$ with $|a|<a_{0}$. It follows that
\[
\sup\limits_{z\in \D_{r}}d_{U}(F\circ f_{n}(t,z),F\circ f_{k}(t,z))\leq \K \sup\limits_{z\in \D_{r}}d_{U}(f_{n}(dt,z),f_{k}(dt,z)).
\]
for all $t\in \s^{1}$. Taking the supremum after $t\in \s^{1}$, we get the desired contraction
\[
d(F(f_{n}), F(f_{k})) \leq \K  d(f_{n},f_{k}),\ \ \ \K<1.
\]
The existence and uniqueness of a fixed point follows from the Banach Fixed Point Theorem.
\qed

The following propositions describe the properties of the fixed point $f^*$.

\begin{prop}\label{prop:holomorphic-z} For any fixed $t\in \s^{1}$, $f^{*}(t,z)=(\varphi_{t}(z),z)$, where $\varphi_{t}:\D_{r}\rightarrow U'$ is holomorphic, and either injective or constant.
\end{prop}
\proof The fixed point $f^{*}$ is obtained via the Banach Fixed Point Theorem as the limit of the sequence $f_{n}(t,z) = F^{\circ n}(f_{0})(t,z)$ for $n\geq 1$ and $f_{0}(t,z)=(\gamma_{0}(t),z)$. We can write $f_{n}(t,z) = (\varphi^{n}_{t}(z),z)$, where $\varphi^{n}_{t}:\D_{r}\rightarrow U'$ are holomorphic and injective for $n\geq 1$. By Hurwitz's theorem a uniform limit of holomorphic injective mappings is holomorphic and either injective or constant.
\qed

\begin{prop}\label{prop:continuous-t}
The function $f^{*}:\s^{1}\times\D_{r}\rightarrow V$ is continuous with respect to $t\in \s^{1}$, holomorphic with respect to $z\in\D_{r}$ and holomorphic with respect to the parameter $a$.
\end{prop}
\proof As observed in the previous proposition, the map $f^{*}$ is obtained as a uniform limit of the sequence $f_{n}(t,z) = (\varphi^{n}_{t}(z),z)$, where
$\varphi^{n}_{t}(z)$ is continuous in $t$ and holomorphic in $z$. Thus $f^{*}$ is continuous in $t$ and holomorphic in $z$. 

Clearly $f_{0}(t,z)=(\gamma_{0}(t),z)$ does not depend on the parameter $a$. When $|a|$ is small,
each function $f_{n}$ depends holomorphically on $a$. The construction of the
metric space is uniform in $a$ and so the limit $f^{*}$ is holomorphic with respect to $a$.
\qed

We can now recover the Julia set $J^{+}\cap V$ as the image of the fixed point $f^{*}$. 
\begin{lemma}\label{lemma:J+}
    $J^{+}\cap V=\bigcap_{n\geq 0}H^{-\circ n}(V)$.
\end{lemma}
\proof Let $q$ be any point in $\bigcap_{n\geq 0}H^{-\circ n}(V)$. Since all forward iterates of $q$ remain in the bounded set $V$, $q$ cannot belong to the escaping set $U^{+}$.  When $H$ is hyperbolic, the interior of $K^{+}$ consists of the basins of attraction of attractive periodic orbits \cite{BS1}. However, the set $U'$ does not contain any attractive cycles of the polynomial $p$ so the set $V=U'\times \D_{r}$ does not contain any attractive cycles of the \He map $H$ for small values of the Jacobian. Since all forward iterates of $q$ remain in $V$, $q$ cannot belong to the interior of $K^{+}$. Hence $q\in J^{+}$.

Let now $q$ be any point in $J^{+}\cap V$. The Julia set $J$ is contained in $V$. When $H$ is hyperbolic, the Julia set $J^+$ is the stable set of $J$, that is $W^{s}(J)=J^{+}$. It follows that $q$ must belong to the stable manifold  $W^{s}(y)$ of some point $y\in J$. So all forward iterates of $q$ converge to the orbit of $y$ which is contained in $J$, hence also in $V$. In particular no forward iterate of $q$ can exit $V$,  hence $q\in \bigcap_{n\geq 0}H^{-\circ n}(V)$.
\qed

\begin{lemma}\label{lemma:image} $Im(f^{*}) = J^{+}\cap V$.
\end{lemma}
\proof 
It is easy to see that $\bigcap_{n\geq 0}H^{-\circ n}(V)=Im(f^*)$ by construction, and that $f^*$ verifies the relation
$H^{-1}(Im(f^{*}))\cap V=Im(f^{*})$. By induction on $n\geq 1$ we get
\begin{equation}\label{eq:im}
H^{-\circ (n+1)}(Im(f^{*}))\cap H^{-\circ n}(V)\cap\ldots\cap H^{-1}(V)\cap V= Im(f^{*}),
\end{equation}
hence $Im(f^{*})\subset \bigcap_{n\geq 0}H^{-\circ n}(V)$.

By Lemma  \ref{lemma:J+} we have
$J^{+}\cap V=\bigcap_{n\geq 0}H^{-\circ n}(V) = \bigcap_{n\geq 0}H^{-\circ n}(V\cap \overline{U^{+}} )$. 
The set $J^+$ is the topological boundary of the set $U^+$ and $J^+ \cap V$ is the inner boundary of the set $V\cap \overline{U^{+}}$.
Recall that $f_{0}(t,z)=(\gamma_{0}(t),z)$.  By construction, $Im(f_{0})$ is the outer boundary of $V$ and is entirely contained in $U^{+}$. 
Moreover, the sequence $f_{n}:\s^{1}\times \D_{r}\rightarrow V$,  $f_n=F^{n}(f_0)$ converges to the fixed point $f^*$. 
The map $f_{n}(\s^{1}\times\D_{r})$ is the outer boundary of the set $\bigcap_{0\leq k\leq n}H^{-\circ k}(V\cap \overline{U^{+}})$. Hence $Im(f^{*})=\bigcap_{n\geq 0}H^{-\circ n}(V\cap \overline{U^{+}})$.
\qed

%================================================================
\section{Characterizations of $J$ and $J^{+}$}\label{sec:conjugacy}

Consider $f^{*}(t,z) = (\varphi_{t}(z),z)$, where $\varphi_{t}(z)$ is continuous with respect to $t\in \s^{1}$ and analytic with respect to $z\in \D_{r}$ and $a$. Let $\sigma:\s^{1}\times \D_{r}\rightarrow \s^{1}\times \D_{r}$  be given by
\begin{equation}\label{eq:sigma}
\sigma(t,z) = \left( dt, a\varphi_{t}(z) \right).
\end{equation}
On the first coordinate this is the $d$-tupling map on the unit circle $t\mapsto dt\ (\mbox{mod} 1)$. 
We chose to disregard the dependency on $a$ in the definition of $\sigma$, to simplify notations.  For sufficiently small $|a|>0$ the map  $\sigma$ is well-defined, open and injective (see \cite[Proposition~12.3]{RT1}). 
Moreover, the map $\sigma$ has the followings expansion with respect to the parameter $a$  (see \cite[Lemma~12.2]{RT1}):
\[
\sigma(t,z) = \left(dt, a\gamma(t) -\frac{a^{2}z}{p'(\gamma(t))}+\bigO(a^{3})\right),
\]
where $\gamma$ is the Carath\'eodory loop of the polynomial $p$.

\begin{thm}\label{thm:conjugacy1} Let $p$ be a hyperbolic
polynomial with connected Julia set.  There exists $a_0>0$ such that
if $0<|a|<a_0$ then the diagram
\[
\diag{\s^{1}\times \D_{r}}{J^{+}\cap V}{\s^{1}\times \D_{r}}{J^{+}\cap V} {f^{*}}{\sigma}{H_{p,a}}{f^{*}}
\]
commutes.
\end{thm}
\proof The existence of the fixed point $f^*$ was established in  Section \ref{sec:spaceF}. By Lemma \ref{lemma:image}, the image of $f^*$ is the set $J^+\cap V$. Consider the definition of $\sigma$ from Equation \ref{eq:sigma}.  We just need to verify the commutativity of the diagram. Since $H\circ f^{*}(t\times \D_{r})$ is compactly contained in $f^{*}(dt\times \D_{r})$ we get that
\[
H\circ f^{*}(t,z)=\hvec{p(\varphi_{t}(z))+az}{a\varphi_{t}(z)} = \hvec{\varphi_{dt}(a\varphi_{t}(z))}{a\varphi_{t}(z)},
\]
which is equal to $f^{*}\circ \sigma (t,z)$ as 
$f^{*}\circ \sigma (t,z) = f^{*}(dt,a\varphi_{t}(z)) = (\varphi_{dt}(a\varphi_{t}(z)),a\varphi_{t}(z))$. 
The last equality holds since $f^{*}(dt\times \D_{r})$ is a vertical-like fiber and can be parametrized by the second coordinate 
via the map $\varphi_{dt}(\cdot)$. 
\qed

Theorem \ref{thm:conjugacy1} gives only a semi-conjugacy between $H$ and $\sigma$, but we are able to identify the equivalence classes of $f^{*}$ explicitly using the fact that $f^*$ is holomorphic with respect to $a$ and $z$ and Hurwitz's theorem (see \cite[Propositions~12.4-12.6]{RT1}):
\begin{equation}\label{eq: equiv-classes}
f^{*}(t_{1},z_1)=f^{*}(t_2,z_2) \mbox{ if and only if }
\gamma(t_1)=\gamma(t_2) \mbox{ and } z_1=z_2.
\end{equation}
This induces a natural equivalence relation on  $\s^{1}\times \D_{r}$: $(t_{1},z)\sim(t_{2},z)$ if and only if $\gamma(t_{1})=\gamma(t_{2})$. 
Notice that in one-dimension this corresponds to the equivalence relation induced by the Thurston lamination on $\s^1$ (see \cite{Th}) which identifies the Julia set $J_p$ to the quotient $\s^1/_{\sim}$.

From relations \ref{eq:sigma} and \ref{eq: equiv-classes} we have that $\varphi_{t_{1}}(z) = \varphi_{t_{1}}(z)$ and  $\sigma(t_{1},z) = \sigma(t_{2},z)$ whenever $\gamma(t_{1})=\gamma(t_{2})$. Therefore, the map $\sigma$ descends to a map on $\s^{1}\times \D_{r}/_{\sim}$. The space $\s^{1}\times \D_{r}/_{\sim}$ is  naturally identified to $J_{p}\times\D_{r}$, so the map $\sigma$ is conjugate to a map $\sigma_{p}$ acting on $J_{p}\times\D_{r}$ of the form 
\[
\sigma_{p}(\zeta,z) = \left(p(\zeta), a\zeta -\frac{a^{2}z}{p'(\zeta)}+\bigO(a^{3})\right).
\]
Note that the map $\sigma_{p}$ is analytic with respect to $z$, which implies that $J^+ \cap V$ is an analytic fiber bundle over $J_p$. We refer to \cite[Section~12]{RT1} for the complete details. 
We can further conjugate $\sigma_{p}$ to a map $\psi: J_{p}\times\D_{r}\rightarrow J_{p}\times\D_{r}$  of the form 
\begin{equation}\label{eq:psi}
\psi (\zeta, z)=\left(p(\zeta),\epsilon\zeta -\frac{\epsilon^{2}z}{p'(\zeta)}\right),
\end{equation}
for some $\epsilon>0$ independent of $a$ (see \cite[Lemmas~12.7, 12.8]{RT1}).

The Julia set $J$ is the set of points from $J^{+}$ that do not escape to infinity under backward iterations of the \He map. By assumption (ii) from the construction of the neighborhood $V$ in Section \ref{sec:V}, $J\subset V$. Thus $J = \bigcap_{n\geq0}H^{\circ n}(J^{+}\cap V)$. 

Let $\Sigma^{+} = \bigcap_{n\geq0}\sigma^{\circ n}(\s^{1}\times \D_{r})$. The following is a direct consequence of Theorem \ref{thm:conjugacy1} and the discussion above: 

\begin{thm}\label{thm:J} The Julia set $J$ of the \He map is homeomorphic to the quotiented solenoid $\Sigma^{+}/_{\sim}$, which is further homeomorphic to $\bigcap_{n\geq0}\psi^{\circ n}(J_{p}\times\D_{r})$. 
\end{thm}

Define the ({\it inductive limit}) space $\widecheck{J_{p}}$ as the quotient of $(J_{p}\times\D_{r})\times \N\big{/}\!\sim$,  where the equivalence relation this time is defined by $(x,n)\sim (\psi(x),n+1)$. The space $\widecheck{J_{p}}$ comes with a natural bijective map $\widecheck{\psi}:\widecheck{J_{p}}\rightarrow \widecheck{J_{p}}$ given by $(x,n)\mapsto (\psi(x),n)$. One should think of the space $\widecheck{J_{p}}$ as an increasing union of sets homeomorphic to $J_{p}\times\D_{r}$. 

\begin{thmHOV}[Hubbard--Oberste-Vorth \cite{HOV2}]\label{thm:HOV} Let $p$ be a  hyperbolic
polynomial with connected Julia set.  There exists $a_0>0$ such that
if $0<|a|<a_0$ then there exists a homeomorphism $\Phi^{+}$
that makes the diagram 
\[
\diag{\widecheck{J_{p}}}{J^{+}}{\widecheck{J_{p}}}{J^{+}} {\Phi^{+}}{\widecheck{\psi}}{H_{p,a}}{\Phi^{+}}
\]
commute.
\end{thmHOV}
\proof We have already shown that there exists a homeomorphism $\Phi^{+}$ that makes the following diagram commute
\begin{equation}\label{eq:diag}
\diag{J_{p}\times \D_{r}}{J^{+}\cap V}{J_{p}\times \D_{r}}{J^{+}\cap V} {\Phi^{+}}{\psi}{H_{p,a}}{\Phi^{+}}
\end{equation}
The map $\Phi^{+}$ is just a composition between the fixed point $f^{*}$ and the homeomorphism that conjugates $\sigma_{p}$ to $\psi$. 
Taking the inductive limit on both sides of the diagram and using \cite[Proposition~6.1]{HOV2} completes the proof. By a small abuse of notation, we use $\Phi^{+}$ to denote the map $\widecheck{\Phi}^{+}$ induced on the inductive limit spaces.
\qed

The map $\psi$ can be conjugate to a map $\psi'$ defined by $\psi'(\zeta, z)=\left(p(\zeta),\zeta -\frac{\epsilon^{2}z}{p'(\zeta)}\right)$ by the linear change of variables $(\zeta,z)\mapsto (\zeta,\epsilon z)$. The map $\psi'$ is in fact the model map used in \cite{HOV2}. The difference comes from the fact that we are using the \He map normalized so that it has Jacobian $-a^{2}$ rather than $a$.

Let $\zeta_{0}\in J_{p}$. Define the ({\it projective limit}) set $\widehat{J}_{p}$ to be the set of pre-histories of a point $\zeta_{0}\in J_{p}$ under the polynomial $p$:
\[
\widehat{J}_{p} = \{ (\zeta_{0},\zeta_{-1},\zeta_{-2},\ldots) : p(\zeta_{-i}) = \zeta_{-i+1}\ \mbox{for all}\ i\geq 1\}
\]
The space $\widehat{J}_{p}$ comes with a natural bijective map $\widehat{p}: \widehat{J}_{p}\rightarrow \widehat{J}_{p}$, given by 
\[
\widehat{p}((\zeta_{0},\zeta_{-1},\zeta_{-2},\ldots)) = (p(\zeta_{0}),p(\zeta_{-1}),p(\zeta_{-2}),\ldots) = (p(\zeta_{0}),\zeta_{0},\zeta_{-1},\ldots).
\]
Let $\widehat{\psi}$ be the map that associates to $(\zeta_{0},\zeta_{-1},\zeta_{-2},\ldots)\in \widehat{J_{p}}$ a point $\zeta^{*}$, the unique point of the intersection $\bigcap_{i\geq 0}\psi^{\circ i}(\zeta_{-i}\times\D_{r})$. Let $\Phi$ denote the composition of $\Phi^{+}$ from diagram \ref{eq:diag} and $\widehat{\psi}$. 

The map $\Phi$ is a homeomorphism from $\widehat{J_{p}}$ to $J$, the Julia set of the \He map, which makes the diagram 
\begin{equation}\label{eq:HOV}
\diag{\widehat{J_{p}}}{J}{\widehat{J_{p}}}{J} {\Phi}{\widehat{p}}{H_{p,a}}{\Phi}
\end{equation}
commute. We have just obtained the same model for the Julia set $J$ as in \cite{HOV2}, which is of course homeomorphic to the model from Theorem \ref{thm:J}. The space $\widehat{J_{p}}$ is a combinatorial model, while the model in Theorem \ref{thm:J} is topological. 

The following theorem gives another perspective on the set $J^+$, without using the inductive limit space $\widecheck{J_p}$ (see also \cite[Theorem~5.2]{R} and \cite{T} for other characterizations).

\begin{thm}\label{thm:conjugacy2} Let $p$ be a hyperbolic
polynomial with connected Julia set.  There exists $a_0>0$ such that
for all parameters $a$ with $0<|a|<a_0$  there exists a semi-conjugacy $\Psi:J_p\times\C\rightarrow J^{+}$ which makes the diagram
\[
\diag{J_{p}\times \C}{J^{+}}{J_{p}\times\C}{J^{+}} {\Psi}{\psi}{H_{p,a}}{\Psi}
\]
commute.
\end{thm}
\proof Let $(\zeta, z)\in J_{p}\times \C$ and let $n$ be the first iterate such that $\psi^{\circ n}(\zeta, z)$ belongs to $J_{p}\times \D_{r}$.
We define $\Psi(\zeta,z) = H_{p,a}^{-\circ n}\circ \Phi^{+}\circ \psi^{\circ n}(\zeta, z)$, where $\Phi^{+}$ is the conjugating homeomorphism from diagram \ref{eq:diag}. It is easy to check that the map $\Psi$ is a surjective semi-conjugacy.
\qed

%=============================================================================
%=============================================================================

%=============================================================================
%=============================================================================
\end{document}